\documentclass[12]{article}
\usepackage{fleqn}
\usepackage{graphicx}
\begin{document}
\Large
\begin{center}{\bf
On the Fine Structure of the Projective Line Over\\ GF(2)$\otimes$GF(2)$\otimes$GF(2)
}
\end{center}
\vspace*{-.0cm}
\begin{center}
Metod Saniga$^{\dag \ddag}$ and Michel Planat$^{\ddag}$
\end{center}
\vspace*{.1cm} \normalsize
\begin{center}
$^{\dag}$Astronomical Institute, Slovak Academy of Sciences\\
SK-05960 Tatransk\' a Lomnica, Slovak Republic\\
(msaniga@astro.sk)

\vspace*{.1cm}
 and

\vspace*{.1cm} $^{\ddag}$Institut FEMTO-ST, CNRS, D\' epartement LPMO, 32 Avenue de
l'Observatoire\\ F-25044 Besan\c con, France\\
(planat@lpmo.edu)
\end{center}

\vspace*{-.2cm} \noindent \hrulefill

\vspace*{.1cm} \noindent {\bf Abstract}

\noindent
The paper gives a succinct appraisal of the properties of the projective
line defined over the direct product ring $R_{\triangle} \equiv$  GF(2)$\otimes$GF(2)$\otimes$GF(2).
The ring is remarkable in that except for unity, all the remaining seven elements are zero-divisors, the non-trivial ones
forming two distinct sets of three; elementary (`slim') and composite (`fat'). Due to this fact, the line in question
is endowed with a very intricate structure. It contains
twenty-seven points, every point has eighteen neighbour points, the neighbourhoods of two distant points share
twelve points and those of three pairwise distant points have six points in common. Algebraically, the points of the 
line can be partitioned into three groups:
a) the group comprising three distinguished points of the ordinary projective line of order two (the `nucleus'), b)
the group composed of twelve points whose coordinates feature both the unit(y) and a zero-divisor (the `inner shell') and c)
the group of twelve points whose coordinates have both the entries zero-divisors (the `outer shell'). The points of the last
two groups can further be split into two subgroups of six points each; while in the former case there is a perfect symmetry
between the two subsets, in the latter case the subgroups have a different footing, reflecting the existence of the two kinds of a zero-divisor.
The structure of the two shells, the way how they are interconnected and their link with the nucleus are all
fully revealed and illustrated in terms of the neighbour/distant relation. Possible applications of this finite ring geometry are also mentioned.
\\




\vspace*{-.1cm} \noindent \hrulefill

\vspace*{.3cm}  \noindent
\section{Introduction}
Projective spaces defined over {\it rings} \cite{tv}--\cite{hav}
represent a very important branch of algebraic geometry, being
once endowed with a richer and more intricate structure when
compared to that of the corresponding spaces defined over {\it
fields}. This difference is perhaps most pronounced in the case of
one-dimensional projective spaces, i.\,e. lines, where those
defined over fields are said to have virtually no intrinsic
structure (see, e.\,g., \cite{bh1}). The projective line defined
over the ring $R_{\triangle} \equiv$
GF(2)$\otimes$GF(2)$\otimes$GF(2), and hereafter denoted as
PR$_{\triangle}$(1), serves as an especially nice illustration of
this point, visible already at the level of cardinalities and
fully revealed when it comes to the neighbour and/or distant
relation. As per the former case, although the underlying ring has
only eight elements, the line itself possesses as many as
twenty-seven points in total, which is to be compared with only
nine points of the line defined over the eight-element field,
GF(8) (see, e.\,g., \cite{de}). Concerning the latter issue,
whilst the neighbour relation is a mere identity relation for any
field line, it acquires, as we shall demonstrate in detail, a
highly non-trivial character for PR$_{\triangle}$(1). The paper,
when combined with \cite{sploc1} and \cite{spdn}, may serve as an
elementary and self-contained introduction into the theory of
finite projective ring geometries aimed mainly at physicists and
scholars of natural sciences.

\section{Rudiments of Ring Theory}
In this section we recollect some basic definitions, concepts and
properties of rings (see, e.\,g., \cite{fr}--\cite{ra}) to be
employed in the sequel.

A {\it ring} is a set $R$ (or, more specifically, ($R,+,*$)) with
two binary operations, usually called addition ($+$) and
multiplication ($*$), such that $R$ is an abelian group under
addition and a semigroup under multiplication, with multiplication
being both left and right distributive over addition.\footnote{It
is customary to denote multiplication in  a ring simply by
juxtaposition, using $ab$ in place of $a*b$.} A ring in which the
multiplication is commutative is a commutative ring. A ring $R$
with a multiplicative identity 1 such that 1$r$ = $r$1 = $r$ for all $r
\in R$ is a ring with unity. A ring containing a finite number of
elements is a finite ring. In what follows the word ring will
always mean a commutative ring with unity.

An element $r$ of the ring $R$ is a {\it unit} (or an invertible
element) if there exists an element $r^{-1}$ such that $rr^{-1} =
r^{-1} r=1$. This element, uniquely determined by $r$,  is called
the multiplicative inverse of $r$. The set of units forms a group
under multiplication. A (non-zero) element $r$ of $R$ is said to
be a (non-trivial) {\it zero-divisor} if there exists $s \neq 0$
such that $sr= rs=0$. An element of a finite ring is either a unit
or a  zero-divisor. A ring in which every non-zero element is a
unit is a {\it field}; finite (or Galois) fields, often denoted by
GF($q$), have $q$ elements and exist only for $q = p^{n}$, where
$p$ is a prime number and $n$ a positive integer. The smallest
positive integer $s$ such that $s1=0$, where $s1$ stands for $1 +
1 + 1 + \ldots + 1$ ($s$ times), is called the {\it
characteristic} of $R$; if $s1$ is never zero, $R$ is said to be
of characteristic zero.

An {\it ideal} ${\cal I}$ of $R$ is a subgroup of $(R,+)$ such
that $a{\cal I} = {\cal I}a \subseteq {\cal I}$ for all $a \in R$.
An ideal of the ring $R$ which is not contained in any other ideal
but $R$ itself is called a {\it maximal}  ideal. If an ideal is of
the form $Ra$ for some element $a$ of $R$ it is called a {\it
principal} ideal, usually denoted by $\langle a \rangle$. A ring
with a unique maximal ideal is a {\it local} ring. Let $R$ be a
ring and ${\cal I}$ one of its ideals. Then $\overline{R} \equiv
R/{\cal I} = \{a + {\cal I} ~|~ a \in R\}$ together with addition
$(a + {\cal I}) + (b + {\cal I}) = a + b +  {\cal I}$ and
multiplication $(a + {\cal I})(b + {\cal I}) = ab +  {\cal I}$ is
a ring, called the quotient, or factor, ring of $R$ with respect
to ${\cal I}$; if ${\cal I}$ is maximal, then $\overline{R}$ is a
field. A very important ideal of a ring is that represented by the
intersection of all maximal ideals; this ideal is called the {\it
Jacobson radical}.

A mapping $\pi$:~ $R \mapsto S$ between two rings $(R,+,*)$ and
$(S,\oplus, \otimes)$ is a ring {\it homomorphism} if it meets the
following constraints: $\pi (a + b) = \pi (a) \oplus \pi (b)$,
$\pi (a * b) = \pi (a) \otimes \pi(b)$ and $\pi (1) = 1$  for any
two elements $a$ and $b$ of $R$. From this definition it is
readily discerned  that $\pi(0) = 0$, $\pi(-a) = -\pi(a)$, a unit
of $R$ is sent into a unit of $S$ and the set of elements $\{a \in
R~ |~ \pi(a) = 0\}$, called the {\it kernel} of $\pi$, is an ideal
of $R$. A {\it canonical}, or {\it natural}, map
$\overline{\pi}$:~$R \rightarrow \overline{R} \equiv R/{\cal I}$
defined by $\overline{\pi}(r) = r + {\cal I}$ is clearly a ring
homomorphism with kernel ${\cal I}$. A bijective ring homomorphism
is called a ring {\it iso}morphism; two rings $R$ and $S$ are
called isomorphic, denoted by $R \cong S$, if there exists a ring
isomorphism between them.

Finally, we mention a couple of relevant examples of rings: a
polynomial ring, $R[x]$, viz. the set of all polynomials in one
variable $x$ and with coefficients in a ring $R$, and the ring
$R_{\otimes}$ that is a (finite) direct product of rings,
$R_{\otimes} \equiv R_{1} \otimes R_{2} \otimes \ldots \otimes
R_{n}$, where both addition and multiplication are carried out
componentwise and where the individual rings need not be the same.

\section{The Ring $R_{\triangle}$, its Fundamental Quotient Rings and\\ Canonical Homomorphisms}
The ring $R_{\triangle} \equiv$ GF(2)$\otimes$GF(2)$\otimes$GF(2)
is, like GF(2) itself, of characteristic two and consists of the
following eight elements
\begin{eqnarray}
R_{\triangle} &=& \{ [0, 0, 0] \equiv a, [1, 0, 0] \equiv b, [0,
1, 0] \equiv c, [0, 0, 1] \equiv d, \nonumber \\ &&
\hspace*{2cm}[1, 1, 0] \equiv e, [1, 0, 1] \equiv f,  [0, 1, 1]
\equiv g, [1, 1, 1] \equiv h \}
\end{eqnarray}
which comprise just one unit,
\begin{equation}
R_{\triangle}^{*} = \{ h = b + c + d \},
\end{equation}
and seven zero-divisors,
\begin{equation}
R_{\triangle} \backslash R_{\triangle}^{*}  = \{a, b, c, d, e = b
+ c , f = b + d, g = c + d \}.
\end{equation}
The addition and multiplication read, respectively, as follows
\begin{center}
\begin{tabular}{||l|cccccccc||}
\hline \hline
$\oplus$ & $a$ & $b$ & $c$ & $d$ & $e$ & $f$ & $g$ & $h$ \\
\hline
$a$ & $a$ & $b$ & $c$ & $d$ & $e$ & $f$ & $g$ & $h$  \\
$b$ & $b$ & $a$ & $e$ & $f$ & $c$ & $d$ & $h$ & $g$ \\
$c$ & $c$ & $e$ & $a$ & $g$ & $b$ & $h$ & $d$ & $f$ \\
$d$ & $d$ & $f$ & $g$ & $a$ & $h$ & $b$ & $c$ & $e$ \\
$e$ & $e$ & $c$ & $b$ & $h$ & $a$ & $g$ & $f$ & $d$ \\
$f$ & $f$ & $d$ & $h$ & $b$ & $g$ & $a$ & $e$ & $c$ \\
$g$ & $g$ & $h$ & $d$ & $c$ & $f$ & $e$ & $a$ & $b$ \\
$h$ & $h$ & $g$ & $f$ & $e$ & $d$ & $c$ & $b$ & $a$ \\
\hline \hline
\end{tabular}~~~~~~
\begin{tabular}{||l|cccccccc||}
\hline \hline
$\otimes$ & $a$ & $b$ & $c$ & $d$ & $e$ & $f$ & $g$ & $h$ \\
\hline
$a$ & $a$ & $a$ & $a$ & $a$ & $a$ & $a$ & $a$ & $a$  \\
$b$ & $a$ & $b$ & $a$ & $a$ & $b$ & $b$ & $a$ & $b$ \\
$c$ & $a$ & $a$ & $c$ & $a$ & $c$ & $a$ & $c$ & $c$ \\
$d$ & $a$ & $a$ & $a$ & $d$ & $a$ & $d$ & $d$ & $d$ \\
$e$ & $a$ & $b$ & $c$ & $a$ & $e$ & $b$ & $c$ & $e$ \\
$f$ & $a$ & $b$ & $a$ & $d$ & $b$ & $f$ & $d$ & $f$ \\
$g$ & $a$ & $a$ & $c$ & $d$ & $c$ & $d$ & $g$ & $g$ \\
$h$ & $a$ & $b$ & $c$ & $d$ & $e$ & $f$ & $g$ & $h$ \\
\hline \hline
\end{tabular}
\end{center}
from where we readily discern that $a \equiv 0$ and $h \equiv 1$
and find out that the ring has three maximal --- and principal as
well --- ideals
\begin{equation}
{\cal I}_{e} \equiv \langle e \rangle = \{ a, b, c, e \},
\end{equation}
\begin{equation}
{\cal I}_{f} \equiv \langle f \rangle = \{ a, b, d, f \},
\end{equation}
and
\begin{equation}
{\cal I}_{g} \equiv \langle g \rangle = \{ a, c, d, g \},
\end{equation}
and three other principal ideals
\begin{equation}
{\cal I}_{1} \equiv \langle b \rangle = \{ a, b \} = {\cal I}_{e}
\cap {\cal I}_{f},
\end{equation}
\begin{equation}
{\cal I}_{2} \equiv \langle c \rangle = \{ a, c \} = {\cal I}_{e}
\cap {\cal I}_{g},
\end{equation}
and
\begin{equation}
{\cal I}_{3} \equiv \langle d \rangle = \{ a, d \} = {\cal I}_{f}
\cap {\cal I}_{g}.
\end{equation}
These ideals give rise to the fundamental quotient rings, all of
characteristic two, namely
\begin{equation}
R_{\triangle}/{\cal I}_{e} \cong R_{\triangle}/{\cal I}_{f} \cong
R_{\triangle}/{\cal I}_{g} \cong {\rm GF(2)}
\end{equation}
and
\begin{equation}
R_{\triangle}/{\cal I}_{1} \cong R_{\triangle}/{\cal I}_{2} \cong
R_{\triangle}/{\cal I}_{3} \cong {\rm GF(2)} \otimes {\rm GF(2)},
\end{equation}
which yield two canonical homomorphisms
\begin{equation}
R_{\triangle} \rightarrow {\rm GF(2)}
\end{equation}
and
\begin{equation}
R_{\triangle} \rightarrow {\rm GF(2)} \otimes {\rm GF(2)},
\end{equation}
respectively. To conclude this section, it is worth mentioning
that there exist two kinds of a subring of $R_{\triangle}$ which
are isomorphic to GF(2)$\otimes$GF(2), differing from each other
in whether or not is their unity inherited from  $R_{\triangle}$.
As per the former case, an example is furnished by the subset
$R_{\circ} \equiv \{ a, b, g, h \}$ with the addition and
multiplication inherited from the parent ring, viz.
\begin{center}
\begin{tabular}{||l|cccc||}
\hline \hline
$\oplus$ & $a$ & $b$ & $g$ & $h$ \\
\hline
$a$ & $a$ & $b$ & $g$ & $h$  \\
$b$ & $b$ & $a$ & $h$ & $g$  \\
$g$ & $g$ & $h$ & $a$ & $b$ \\
$h$ & $h$ & $g$ & $b$ & $a$ \\
\hline \hline
\end{tabular}~~~~~~~
\begin{tabular}{||l|cccc||}
\hline \hline
$\otimes$ & $a$ & $b$ & $g$ & $h$  \\
\hline
$a$ & $a$ & $a$ & $a$ & $a$  \\
$b$ & $a$ & $b$ & $a$ & $b$  \\
$g$ & $a$ & $a$ & $g$ & $g$  \\
$h$ & $a$ & $b$ & $g$ & $h$  \\
\hline \hline
\end{tabular}~.
\end{center}
That $R_{\circ} \cong$ GF(2)$\otimes$GF(2) stems from the
following correspondence $a =[0, 0]$, $b =[1, 0]$, $g =[0, 1]$ and
$h =[1, 1]$, and we see that $h$ is indeed the unity in both
$R_{\triangle}$ and $R_{\circ}$. An example of the latter case is
provided by the subset $R_{\bullet} \equiv \{ a, b, c, e\}$, with
the following addition and multiplication tables
\begin{center}
\begin{tabular}{||l|cccc||}
\hline \hline
$\oplus$ & $a$ & $b$ & $c$ & $e$ \\
\hline
$a$ & $a$ & $b$ & $c$ & $e$  \\
$b$ & $b$ & $a$ & $e$ & $c$  \\
$c$ & $c$ & $e$ & $a$ & $b$ \\
$e$ & $e$ & $c$ & $b$ & $a$ \\
\hline \hline
\end{tabular}~~~~~~~
\begin{tabular}{||l|cccc||}
\hline \hline
$\otimes$ & $a$ & $b$ & $c$ & $e$  \\
\hline
$a$ & $a$ & $a$ & $a$ & $a$  \\
$b$ & $a$ & $b$ & $a$ & $b$  \\
$c$ & $a$ & $a$ & $c$ & $c$  \\
$e$ & $a$ & $b$ & $c$ & $e$  \\
\hline \hline
\end{tabular}~.
\end{center}
We can easily verify that $R_{\bullet} \cong$ GF(2)$\otimes$GF(2)
by making the following identifications $a =[0, 0]$, $b =[1, 0]$,
$c =[0, 1]$ and $e =[1, 1]$, and notice that the current unity,
$e$, is  a {\it zero-divisor} in $R_{\triangle}$.

\section{The Projective Line over $R_{\triangle}$ and its Fine Structure}
Given a ring $R$ and GL$_{2}$($R$), the general linear group of
invertible two-by-two matrices with entries in $R$, a pair
($\alpha, \beta$) $\in R^{2}$ is called {\it admissible} over $R$
if there exist $\gamma, \delta \in R$ such that \cite{her}
\begin{equation}
\left(
\begin{array}{cc}
\alpha & \beta \\
\gamma & \delta \\
\end{array}
\right) \in {\rm GL}_{2}(R).
\end{equation}
The projective line over $R$, PR(1), is defined as the set of
classes of ordered pairs $(\varrho \alpha, \varrho \beta)$, where
$\varrho$ is a unit and $(\alpha, \beta)$ admissible
\cite{bh1},\cite{hav},\cite{her},\cite{bh2}. Such a line carries
two non-trivial, mutually complementary relations of neighbour and
distant. In particular, its two distinct points $X$: $(\varrho
\alpha, \varrho \beta)$ and $Y$: $(\varrho \gamma, \varrho
\delta)$ are called {\it neighbour} (or, also {\it parallel}) if
\begin{equation}
\left(
\begin{array}{cc}
\alpha & \beta \\
\gamma & \delta \\
\end{array}
\right) \notin {\rm GL}_{2}(R)
\end{equation}
and {\it distant} otherwise, i.\,e., if condition (14) is met. The
neighbour relation is reflexive (every point is obviously
neighbour to itself) and symmetric (i.\,e., if $X$ is neighbour to
$Y$ then also $Y$ is neighbour to $X$), but, in general, not
transitive (i.\,e., $X$  being neighbour to $Y$ and $Y$ being
neighbour to $Z$ does not necessarily mean that $X$ is neighbour
to $Z$ --- see, e.\,g., \cite{vk95},\cite{hav}). Given a point of
PR(1), the set of all neighbour points to it will be called its
{\it neighbourhood}.\footnote{To avoid any confusion, the reader
must be warned here that some authors (e.\,g.
\cite{hav},\cite{bh2}) use this term for the set of {\it distant}
points instead.} Obviously, if $R£$ is a field then `neighbour'
simply reduces to `identical' and `distant' to `different'.

Our next task is to apply this general definitions and concepts to
our ring $R_{\triangle}$. To this end in view, one first notice
that Eqs.\,(14) and (15) acquire, respectively, the following
simple forms
\begin{equation}
\det \left(
\begin{array}{cc}
\alpha & \beta \\
\gamma & \delta \\
\end{array}
\right) = 1
\end{equation}
and
\begin{equation}
\det \left(
\begin{array}{cc}
\alpha & \beta \\
\gamma & \delta \\
\end{array}
\right) \in R_{\triangle} \backslash R_{\triangle}^{*}.
\end{equation}
Employing the former of the two equations, we find that
PR$_{\triangle}$(1) consists of the following twenty-seven points:

1) the three distinguished points (the `nucleus')
\begin{equation}
U:~ (1, 0),~~ V:~ (0, 1),~~ W:~ (1, 1),
\end{equation}
which represent the ordinary projective line of order two
(PG(1,2)) embedded in PR$_{\triangle}$(1);

2) the twelve points of the `inner shell' whose coordinates
feature both the unity and a zero-divisor,
\begin{eqnarray}
&&I_{1}^{S}:~(1,b),~~I_{2}^{S}:~(1,c),~~I_{3}^{S}:~(1,d),~~I_{1}^{F}:~(1,e),~~I_{2}^{F}:~(1,f),~~I_{3}^{F}:~(1,g),
\nonumber \\
&& \\
&&J_{1}^{S}:~(b,1),~~J_{2}^{S}:~(c,1),~~J_{3}^{S}:~(d,1),~~J_{1}^{F}:~(e,1),~~J_{2}^{F}:~(f,1),~~J_{3}^{F}:~(g,1),
\nonumber
\end{eqnarray}
and which form two symmetric sets of six points each, the sets
themselves subject to further splitting according as the
zero-divisor is elementary (slim, `S') or composite (fat, `F') ---
see Eq.\,(3); and

3) the twelve points of the `outer shell' whose coordinates have
zero-divisors in both the entries,
\begin{eqnarray}
&&S_{1}^{+}:~(e,d),~~S_{2}^{+}:~(f,c),~~S_{3}^{+}:~(g,b),~~S_{1}^{-}:~(d,e),~~S_{2}^{-}:~(c,f),~~S_{3}^{-}:~(b,g),
\nonumber \\
&& \\
&&F_{1}^{+}:~(e,f),~~F_{2}^{+}:~(e,g),~~F_{3}^{+}:~(f,g),~~F_{1}^{-}:~(f,e),~~F_{2}^{-}:~(g,e),~~F_{3}^{-}:~(g,f),
\nonumber
\end{eqnarray}
and which comprise two asymmetric sets of cardinality six
according as both the entries are composite zero-divisors (`F') or
not (`S'); here the sets are further subdivided in terms of the
parity (`$+$' or `$-$') of the coordinates' entries. With the help
of Eq.\,(17) we can discover how these three sets, and elements
within them, are related to each other.

To get the idea about the general structure of the line, we first
consider the nucleus and find out that these three mutually
distant points have the following neighbourhoods ($i=1, 2, 3$):
\begin{equation}
U:~~\{I_{i}^{S},~I_{i}^{F},~S_{i}^{+},~S_{i}^{-},~F_{i}^{+},~F_{i}^{-}\},
\end{equation}
\begin{equation}
V:~~\{J_{i}^{S},~J_{i}^{F},~S_{i}^{+},~S_{i}^{-},~F_{i}^{+},~F_{i}^{-}\},
\end{equation}
\begin{equation}
W:~~\{I_{i}^{S},~I_{i}^{F},~J_{i}^{S},~J_{i}^{F},~F_{i}^{+},~F_{i}^{-}\}.
\end{equation}
Now, as the coordinate system on this line can {\it always} be
chosen in such a way that the coordinates of {\it any} three
mutually distant points are made identical to those of $U$, $V$
and $W$, from the last three expressions we discern that the
neighbourhood of any  point of the line features eighteen distinct
points, the neighbourhoods of any two distant points share twelve
points and the neighbourhoods of any three mutually distant points
have six points in common --- as depicted in Fig.\,1. As in the
case of the lines defined over GF(2)[$x$]/$\langle x^{3} - x
\rangle$ and GF(2) $\otimes$ GF(2) \cite{sploc1}, the neighbour
relation is not transitive; however, a novel feature, not
encountered in the previous cases, is here a non-zero overlapping
between the neighbourhoods of {\it three} pairwise distant points,
which can be attributed to the existence of three maximal ideals
of $R_{\triangle}$.
\begin{figure}[t]
\centerline{\includegraphics[width=8.6truecm,clip=]{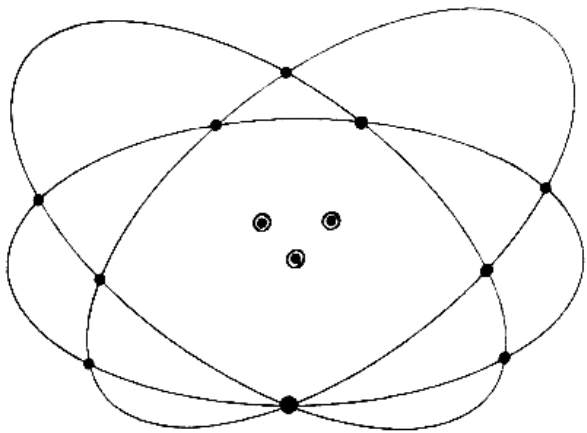}}
\caption{A schematic sketch of the structure of the projective
line $PR_{\triangle}(1)$. Choosing any three pairwise distant
points (represented by the three double circles), the remaining
points of the line are all located on the neighbourhoods of the
three points (three sets of points located on three different
ellipses centered at the points in question). Every small bullet
represents {\it two} distinct points of the line, while the big
bullet at the bottom stands for as many as {\it six} different
points.}
\end{figure}

\noindent The intricacies of the geometry are fully revealed if
one examines the neighbour/distant relation between the points of
the shells as summarized in the following tables (Tables 1--3),
where the sign `+/$-$' means, respectively,  distant/neighbour;
thus, for example, from the first table we read off that the
points $J_{1}^{S}$ and $I_{3}^{F}$ are distant, whilst from the
third table we find out that $F_{2}^{+}$ and $J_{2}^{S}$ are
neighbour. We readily see that the inner shell exhibits a more
pronounced asymmetry in the coupling between the individual
elements than the outer one, as it is also obvious from the
comparison of the two parts of Fig.\,2. Further, as per the
interconnection between the two shells, one finds out that whereas
the `fat' points of the outer shell are symmetrically
(two-and-two) coupled with the points of the two sets of the inner
shell, the `slim' points show a strong (three-and-one) asymmetry
in this respect. Noteworthy is also a very intricate character of
the coupling between the `fat' points of the outer shell, which
sharply contrasts the triviality of the corresponding coupling in
the inner shell. To complete the picture, we also give the table
(Table 4) showing the connection between the nucleus and the
shells. There are obviously a number of other interesting
sub-configurations of the line, like a seven-point configuration
comprising a point of the nucleus and one of the sets of either
shell and a fifteen-point configuration featuring the nucleus and
a shell.
\begin{table}
\caption{The neighbour/distant relation between the points of the
inner shell of the projective line PR$_{\triangle}$(1).}
\begin{center}
\begin{tabular}{||l|ccc|ccc|ccc|ccc||}
\hline \hline
 & $I_{1}^{S}$ & $I_{2}^{S}$ & $I_{3}^{S}$ & $I_{1}^{F}$ & $I_{2}^{F}$ & $I_{3}^{F}$
 & $J_{1}^{S}$ & $J_{2}^{S}$ & $J_{3}^{S}$ & $J_{1}^{F}$ & $J_{2}^{F}$ & $J_{3}^{F}$ \\
\hline
$I_{1}^{S}$ & $-$ & $-$& $-$& $-$& $-$& $+$ & $-$& $+$ & $+$ & $-$& $-$& $+$ \\
$I_{2}^{S}$ & $-$ & $-$& $-$& $-$& $+$& $-$ & $+$& $-$ & $+$ & $-$& $+$& $-$ \\
$I_{3}^{S}$ & $-$ & $-$& $-$& $+$& $-$& $-$ & $+$& $+$ & $-$ & $+$& $-$& $-$ \\
\hline
$I_{1}^{F}$ & $-$ & $-$& $+$& $-$& $-$& $-$ & $-$& $-$ & $+$ & $-$& $-$& $-$ \\
$I_{2}^{F}$ & $-$ & $+$& $-$& $-$& $-$& $-$ & $-$& $+$ & $-$ & $-$& $-$& $-$ \\
$I_{3}^{F}$ & $+$ & $-$& $-$& $-$& $-$& $-$ & $+$& $-$ & $-$ & $-$& $-$& $-$ \\
\hline
$J_{1}^{S}$ & $-$ & $+$& $+$& $-$& $-$& $+$ & $-$& $-$ & $-$ & $-$& $-$& $+$ \\
$J_{2}^{S}$ & $+$ & $-$& $+$& $-$& $+$& $-$ & $-$& $-$ & $-$ & $-$& $+$& $-$ \\
$J_{3}^{S}$ & $+$ & $+$& $-$& $+$& $-$& $-$ & $-$& $-$ & $-$ & $+$& $-$& $-$ \\
\hline
$J_{1}^{F}$ & $-$ & $-$& $+$& $-$& $-$& $-$ & $-$& $-$ & $+$ & $-$& $-$& $-$ \\
$J_{2}^{F}$ & $-$ & $+$& $-$& $-$& $-$& $-$ & $-$& $+$ & $-$ & $-$& $-$& $-$ \\
$J_{3}^{F}$ & $+$ & $-$& $-$& $-$& $-$& $-$ & $+$& $-$ & $-$ & $-$& $-$& $-$ \\
\hline \hline
\end{tabular}
\end{center}
\end{table}

\vspace*{.5cm}
\begin{table}
\caption{The neighbour/distant relation between the points of the
outer shell of the projective line PR$_{\triangle}$(1).}
\begin{center}
\begin{tabular}{||l|ccc|ccc|ccc|ccc||}
\hline \hline
 & $F_{1}^{+}$ & $F_{2}^{+}$ & $F_{3}^{+}$ & $F_{1}^{-}$ & $F_{2}^{-}$ & $F_{3}^{-}$
 & $S_{1}^{+}$ & $S_{2}^{+}$ & $S_{3}^{+}$ & $S_{1}^{-}$ & $S_{2}^{-}$ & $S_{3}^{-}$ \\
\hline
$F_{1}^{+}$ & $-$ & $-$& $+$& $-$& $+$& $-$ & $-$& $+$ & $-$ & $+$& $-$& $-$ \\
$F_{2}^{+}$ & $-$ & $-$& $-$& $+$& $-$& $+$ & $-$& $-$ & $+$ & $+$& $-$& $-$ \\
$F_{3}^{+}$ & $+$ & $-$& $-$& $-$& $+$& $-$ & $-$& $-$ & $+$ & $-$& $+$& $-$ \\
\hline
$F_{1}^{-}$ & $-$ & $+$& $-$& $-$& $-$& $+$ & $+$& $-$ & $-$ & $-$& $+$& $-$ \\
$F_{2}^{-}$ & $+$ & $-$& $+$& $-$& $-$& $-$ & $+$& $-$ & $-$ & $-$& $-$& $+$ \\
$F_{3}^{-}$ & $-$ & $+$& $-$& $+$& $-$& $-$ & $-$& $+$ & $-$ & $-$& $-$& $+$ \\
\hline
$S_{1}^{+}$ & $-$ & $-$& $-$& $+$& $+$& $-$ & $-$& $-$ & $-$ & $+$& $-$& $-$ \\
$S_{2}^{+}$ & $+$ & $-$& $-$& $-$& $-$& $+$ & $-$& $-$ & $-$ & $-$& $+$& $-$ \\
$S_{3}^{+}$ & $-$ & $+$& $+$& $-$& $-$& $-$ & $-$& $-$ & $-$ & $-$& $-$& $+$ \\
\hline
$S_{1}^{-}$ & $+$ & $+$& $-$& $-$& $-$& $-$ & $+$& $-$ & $-$ & $-$& $-$& $-$ \\
$S_{2}^{-}$ & $-$ & $-$& $+$& $+$& $-$& $-$ & $-$& $+$ & $-$ & $-$& $-$& $-$ \\
$S_{3}^{-}$ & $-$ & $-$& $-$& $-$& $+$& $+$ & $-$& $-$ & $+$ & $-$& $-$& $-$ \\
\hline \hline
\end{tabular}
\end{center}
\end{table}
\begin{table}[h]
\caption{The neighbour/distant relation between the points of the
two shells of the projective line PR$_{\triangle}$(1).}
\begin{center}
\begin{tabular}{||l|ccc|ccc|ccc|ccc||}
\hline \hline
 & $I_{1}^{S}$ & $I_{2}^{S}$ & $I_{3}^{S}$ & $I_{1}^{F}$ & $I_{2}^{F}$ & $I_{3}^{F}$
 & $J_{1}^{S}$ & $J_{2}^{S}$ & $J_{3}^{S}$ & $J_{1}^{F}$ & $J_{2}^{F}$ & $J_{3}^{F}$ \\
\hline
$F_{1}^{+}$ & $-$ & $+$& $-$& $-$& $-$& $+$ & $-$& $-$ & $+$ & $-$& $-$& $+$ \\
$F_{2}^{+}$ & $+$ & $-$& $-$& $-$& $+$& $-$ & $-$& $-$ & $+$ & $-$& $+$& $-$ \\
$F_{3}^{+}$ & $+$ & $-$& $-$& $+$& $-$& $-$ & $-$& $+$ & $-$ & $+$& $-$& $-$ \\
\hline
$F_{1}^{-}$ & $-$ & $-$& $+$& $-$& $-$& $+$ & $-$& $+$ & $-$ & $-$& $-$& $+$ \\
$F_{2}^{-}$ & $-$ & $-$& $+$& $-$& $+$& $-$ & $+$& $-$ & $-$ & $-$& $+$& $-$ \\
$F_{3}^{-}$ & $-$ & $+$& $-$& $+$& $-$& $-$ & $+$& $-$ & $-$ & $+$& $-$& $-$ \\
\hline
$S_{1}^{+}$ & $-$ & $-$& $-$& $+$& $-$& $-$ & $-$& $-$ & $+$ & $-$& $+$& $+$ \\
$S_{2}^{+}$ & $-$ & $-$& $-$& $-$& $+$& $-$ & $-$& $+$ & $-$ & $+$& $-$& $+$ \\
$S_{3}^{+}$ & $-$ & $-$& $-$& $-$& $-$& $+$ & $+$& $-$ & $-$ & $+$& $+$& $-$ \\
\hline
$S_{1}^{-}$ & $-$ & $-$& $+$& $-$& $+$& $+$ & $-$& $-$ & $-$ & $+$& $-$& $-$ \\
$S_{2}^{-}$ & $-$ & $+$& $-$& $+$& $-$& $+$ & $-$& $-$ & $-$ & $-$& $+$& $-$ \\
$S_{3}^{-}$ & $+$ & $-$& $-$& $+$& $+$& $-$ & $-$& $-$ & $-$ & $-$& $-$& $+$ \\
\hline \hline
\end{tabular}
\end{center}
\end{table}
\begin{table}
\caption{The neighbour/distant relation between the points of the
nucleus and those of the inner/outer shell of the projective line
PR$_{\triangle}$(1).}
\begin{center}
\begin{tabular}{||l|ccc|ccc|ccc|ccc||}
\hline \hline
 & $I_{1}^{S}$ & $I_{2}^{S}$ & $I_{3}^{S}$ & $I_{1}^{F}$ & $I_{2}^{F}$ & $I_{3}^{F}$
 & $J_{1}^{S}$ & $J_{2}^{S}$ & $J_{3}^{S}$ & $J_{1}^{F}$ & $J_{2}^{F}$ & $J_{3}^{F}$ \\
\hline
$U$ & $-$ & $-$& $-$& $-$& $-$& $-$ & $+$& $+$ & $+$ & $+$& $+$& $+$ \\
$V$ & $+$ & $+$& $+$& $+$& $+$& $+$ & $-$& $-$ & $-$ & $-$& $-$& $-$ \\
$W$ & $-$ & $-$& $-$& $-$& $-$& $-$ & $-$& $-$ & $-$ & $-$& $-$& $-$ \\
\hline
\end{tabular}

\vspace*{0.1cm}
\begin{tabular}{||l|ccc|ccc|ccc|ccc||}
\hline
 & $F_{1}^{+}$ & $F_{2}^{+}$ & $F_{3}^{+}$ & $F_{1}^{-}$ & $F_{2}^{-}$ & $F_{3}^{-}$
 & $S_{1}^{+}$ & $S_{2}^{+}$ & $S_{3}^{+}$ & $S_{1}^{-}$ & $S_{2}^{-}$ & $S_{3}^{-}$ \\
\hline
$U$ & $-$ & $-$& $-$& $-$& $-$& $-$ & $-$& $-$ & $-$ & $-$& $-$& $-$ \\
$V$ & $-$ & $-$& $-$& $-$& $-$& $-$ & $-$& $-$ & $-$ & $-$& $-$& $-$ \\
$W$ & $-$ & $-$& $-$& $-$& $-$& $-$ & $+$& $+$ & $+$ & $+$& $+$& $+$ \\
\hline \hline
\end{tabular}
\end{center}
\end{table}
\begin{figure}[t]
\centerline{\includegraphics[width=6.3truecm,clip=]{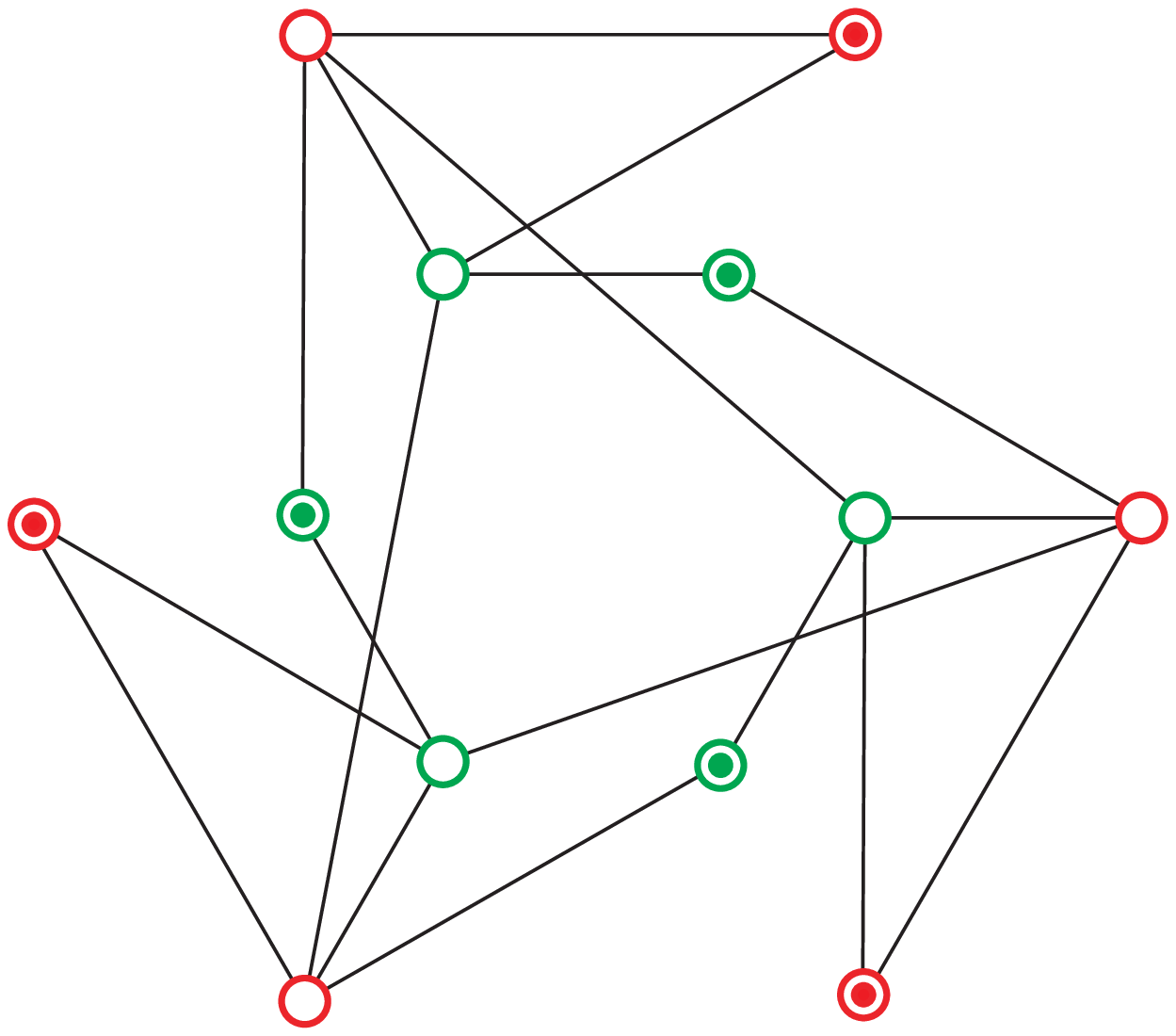}
\includegraphics[width=5truecm,clip=]{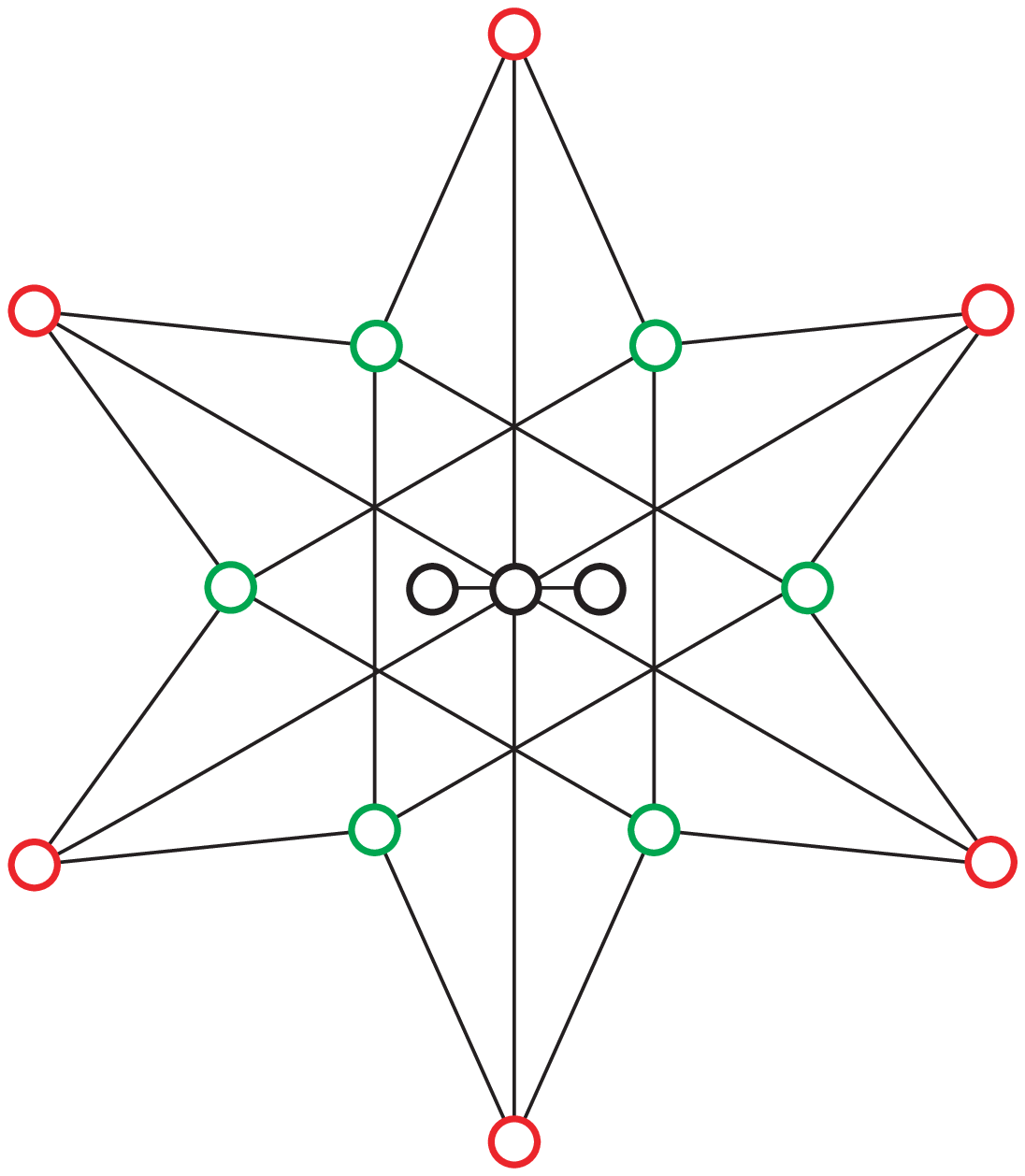}}
\caption{A schematic illustration of the fine structure of the
inner ({\it left}) and outer ({\it right}) shells of the
projective line PR$_{\triangle}$(1), where any two distant points
are joined by a line segment. In both the cases, the two sets of
points are distinguished by different colours; in the former case,
in addition, the filled circles denote the `fat' points and in the
latter case also the nucleus (three points in the center, the
middle one being $W$) and its relation to the shell in question
are shown.}
\end{figure}

\section{Conclusion}
\vspace*{-.4cm} We have carried out an in-depth examination of the
structure of the projective line PR$_{\triangle}$(1). The line was
shown to be endowed, from the algebraic point of view, with three
distinct `layers'; the nucleus comprising the ordinary projective
(sub)line of order two, the inner shell made of points whose
coordinates feature both the unity and a zero-divisor and the
outer shell consisting of the points whose coordinates have
zero-divisors in both the entries. The intricacies of the
structure of the line were fully revealed by employing the
concepts of neighbour/distant, as summarized in Tables 1 to 4 and
partially illustrated in Fig.\,2. We believe that this remarkable
finite ring geometry --- and its multifaceted sub-configurations
--- will soon find interesting applications in quantum physics (to
address, for example, the properties of finite-dimensional quantum
systems following and properly generalizing the strategy pursued
in \cite{sploc2}), chemistry (when dealing with highly-complex
systems of molecules), particle physics, crystallography, biology
and other natural sciences as well.

\vspace*{.6cm} \noindent \Large {\bf Acknowledgements} \normalsize

\vspace*{.0cm} \noindent The first author thanks Mr. P. Bend\'
{\i}k for a careful drawing of the first figure and is extremely
grateful to Dr. R. Kom\v z\' {\i}k for a computer-related
assistance. This work was partially supported by the Science and
Technology Assistance Agency under the contract $\#$
APVT--51--012704, the VEGA project $\#$ 2/6070/26 (both from
Slovak Republic) and by the trans-national ECO-NET project $\#$
12651NJ ``Geometries Over Finite Rings and the Properties of
Mutually Unbiased Bases" (France).


\vspace*{-.1cm}
\begin{thebibliography}{10}
\vspace*{-.20cm}
\bibitem{tv}
T\" orner G, Veldkamp FD. Literature on geometry over rings. J Geom 1991;42:180--200.
\vspace*{-.25cm}
\bibitem{vk95}
Veldkamp FD. Geometry over rings. In: Buekenhout F, editor. Handbook of incidence geometry.
        Amsterdam: Elsevier; 1995. p.\,1033--84.
\vspace*{-.25cm}
\bibitem{vk81}
Veldkamp FD. Projective planes over rings of stable rank 2. Geom Dedicata 1981;11:285--308.
\vspace*{-.25cm}
\bibitem{vk85}
Veldkamp FD. Projective ring planes and their homomorphisms. In: Kaya R, Plaumann P,  Strambach K, editors. Rings and
        geometry (NATO ASI). Dordrecht: Reidel; 1985. p.\,289--350.
\vspace*{-.25cm}
\bibitem{vk84}
Veldkamp FD. Projective ring planes: some special cases. Rend Sem Mat Brescia 1984;7:609--15.
\vspace*{-.25cm}
\bibitem{bh1}
Blunck A, Havlicek H. Projective representations I: Projective lines over a ring.  Abh Math Sem Univ Hamburg 2000;70:287--99.
\vspace*{-.25cm}
\bibitem{hav}
Havlicek H. Divisible designs, Laguerre geometry, and beyond. A preprint available from
$\langle$http://www.geometrie.tuwien.ac.at/havlicek/dd-laguerre.pdf$\rangle$.
\vspace*{-.25cm}
\bibitem{de}
Dembowski P. Finite geometries. Berlin -- New York: Springer;
1968.
\vspace*{-.25cm}
\bibitem{sploc1}
Saniga M, Planat M. The projective line over the finite quotient
ring GF(2)[$x$]/$\langle x^3-x \rangle$ and quantum entanglement
I. Theoretical background. Available from
$\langle$quant-ph/0603051$\rangle$. \vspace*{-.25cm}
\bibitem{spdn}
Saniga M,  Planat M. Projective planes over ``Galois" double
numbers and a geometrical principle of complementarity. J Phys A:
Math Gen 2006, submitted. Available from
$\langle$math.NT/0601261$\rangle$. \vspace*{-.25cm}
\bibitem{fr}
Fraleigh JB. A first course in abstract algebra (5th edition).
Reading (MA): Addison-Wesley; 1994. p.\,273--362. \vspace*{-.25cm}
\bibitem{mcd}
McDonald BR. Finite rings with identity. New York: Marcel Dekker;
1974. \vspace*{-.25cm}
\bibitem{ra}
Raghavendran R. Finite associative rings. Comp Mathematica
1969;21:195--229.
\vspace*{-.25cm}
\bibitem{her}
Herzer A. Chain geometries. In: Buekenhout F, editor. Handbook of
incidence geometry.
        Amsterdam: Elsevier; 1995. p.\,781--842.
\vspace*{-.25cm}
\bibitem{bh2}
Blunck A, Havlicek H. Radical parallelism on projective lines and
non-linear models of affine spaces. Mathematica Pannonica
2003;14:113--27. \vspace*{-.25cm}
\bibitem{sploc2}
Saniga M, Planat M, Minarovjech M. The projective line over the
finite quotient ring GF(2)[$x$]/$\langle x^3-x \rangle$ and
quantum entanglement II. The Mermin ``Magic'' Square/Pentagram.
Available from $\langle$quant-ph/0603206$\rangle$.
\end{thebibliography}
\end{document}